\documentclass[12pt]{article}
\usepackage{cite,amsfonts,amssymb}
\newtheorem{prop}{Proposition}

\newtheorem{theor}{Theorem}

\newtheorem{cor}{Corollary}

\newtheorem{rem}{Remark}

\newtheorem{defin}{Definition}

\newcommand{\bgproof}{\noindent {\bf Proof} \hspace{2mm}}
\newcommand{\edproof}{\hfill $\blacksquare$ \vspace{3mm}}

\title{A series approach to stochastic Volterra equations of convolution type}

\author{{\large\sf Bartosz Bandrowski and Anna Karczewska }\\[2mm]
  \normalsize\it
 Faculty of Mathematics, Computer Science and Econometrics\\ \normalsize\it
 University of Zielona G\'ora\\ \normalsize\it
 ul. Szafrana 4a, 65-246 Zielona G\'ora, Poland\\ \normalsize\it
 e-mail: A.Karczewska@wmie.uz.zgora.pl\\[2mm]
 }

\begin{document}

\maketitle

\def\thefootnote{}
\footnotetext{\noindent {\em 2010 Mathematics Subject
Classification:}
primary: 60H20; secondary: 60H05, 45D05.\\
{\em Key words and phrases:} stochastic linear Volterra equation,
strong solution, resolvent, mild solution, stochastic convolution, series approach.}


\begin{abstract}
In the paper stochastic Volterra equations with noise terms driven by series of 
independent scalar Wiener processes are considered. In our study we use the resolvent approach to the equations under consideration. We give sufficient condition for the existence of strong solution to the class of stochastic Volterra equations of convolution type. We provide regularity of stochastic convolution, as well.
\end{abstract}

\maketitle

\section{Introduction}\label{sSW1}

Let $(H,|\cdot|_H)$ be a seperable Hilbert space and let $(\Omega,\mathcal{F},(\mathcal{F})_{t\ge 0},P)$ denote a probability space.
We consider stochastic Volterra equations in $H$ of the form 
\begin{equation} \label{eSW1}
X(t) = X_0 + \int_0^t a(t-\tau)\,AX(\tau)d\tau +
 \sum\limits_{i=1}^{\infty}\int\limits_0^t \Psi_i(\tau)dW_i(\tau)\;,
\quad t\geq 0\;,
\end{equation}
where $X_0$ is an $H$-valued $\mathcal{F}_0$-measurable random variable, a kernel function 
$a\in L^1_{\mathrm{loc}}(\mathbb{R}_+;\mathbb{R})$ and $A$ is a closed
unbounded linear operator in $H$ with a dense domain $D(A)$. The
domain $D(A)$ is equipped with the graph norm $|\cdot |_{D(A)}$ of
$A$, i.e.\ $|h|_{D(A)}:=(|h|_H^2+|Ah|_H^2)^{1/2}$.

In our work the equation (\ref{eSW1}) is driven by  series of scalar 
Wiener processes;  $W= (W_i)_{i=1}^{\infty}$     and  $\Psi=(\Psi_i)_{i=1}^{\infty} $    are appropriate processes defined below. 

The goal of this paper is to formulate sufficient conditions for the existence and regularity of strong solutions to the Volterra equation driven by infinite dimensionsl Wiener processes. Previously, in \cite{DaPrZa,Ka1,Ka2,KaLi2}, the stochastic integral for Hilbert-Schmidt operator-valued integrands and Wiener processes with values in Hilbert space has been constructed. Moreover, the particular series expansion of the Wiener process with respect to the eigenvectors of its covariance operator has been used. The stochastic integral used in this paper, originally introduced in \cite{Ga}, bases on the construction directly in terms of the sequence of independent scalar processes. In consequence, the stochastic integral is independent of any covariance operator usually connected with a noise process.

The paper is organized as follows. In section \ref{sSW1} we recall the 
deterministic notions and auxiliary results. Section \ref{sSW2} consists of a construction of the stochastic 
integral due to O. van Gaans   \cite{Ga}. In section \ref{sSW3} 
we compare mild and weak solutions and then we provide sufficient condition for stochastic convolution to be a strong solution to the equation (\ref{eSW1}). Section \ref{sSW4} gives regularity of stochastic convolution arising in Volterra equation.
Section \ref{sSW5} contains proof of the main result of the paper formulated in section \ref{sSW3}.

As we have already written, we use the resolvent approach to the equation (\ref{eSW1}).
This means that a deterministic counterpart 
of the equation (\ref{eSW1}), that is, the
equation
\begin{equation} \label{eSW1d}
u(t) = \int_0^t a(t-\tau)\,Au(\tau)d\tau +f(t), \quad t\geq 0,
\end{equation}
admits a resolvent family.
 In (\ref{eSW1d}), the operator $A$ and the
kernel function $a$ are the same as previously considered in (\ref{eSW1}) and $f$ is a $H$-valued function.

By $S(t),~t\geq 0$, we shall denote the family of resolvent operators
corresponding to the Volterra equation (\ref{eSW1d}), 
which is  defined as follows.

\begin{defin}\label{dSW1} (see, e.g.\ \cite{Pr2})\\
A family $(S(t))_{t\geq 0}$ of bounded linear operators in $H$ is
called {\tt resolvent} for (\ref{eSW1d}) if the following
conditions are satisfied:
\begin{enumerate}
\item $S(t)$ is strongly continuous on $\mathbb{R}_+$ and $S(0)=I$;
\item $S(t)$ commutes with the operator $A$, that is, $S(t)(D(A))\subset
D(A)$ and $AS(t)x=S(t)Ax$ for all $x\in D(A)$ and $t\geq 0$;
\item the following {\tt resolvent equation} holds
\begin{equation} \label{eSW2}
S(t)x = x + \int_0^t a(t-\tau) AS(\tau)x d\tau
\end{equation}
for all $x\in D(A),~t\geq 0$.
\end{enumerate}
\end{defin}

Necessary and sufficient conditions for the existence of the
resolvent family have been studied in \cite{Pr2}. This is worth to note 
that the resolvent $S(t),~t\geq 0$, is determined by the operator
$A$ and the function $a$, so we also may say that the pair $(A,a)$
admits a resolvent family.

In the paper we shall assume that the kernel function is completely positive in the following sense.

\begin{defin} \label{dSW2}
We say that function $a\in L^1(0,T)$ is {\tt completely positive}
on $[0,T]$ if for any $\mu\geq 0$, the solutions of the
convolution equations
\begin{equation}
s(t) + \mu \!\int_0^t \! a(t-\tau)s(\tau)d\tau =1 \quad \mbox{and} \quad r(t) + \!\int_0^t \! a(t-\tau)r(\tau)d\tau= a(t)
\end{equation}
are nonnegative on $[0,T]$.
\end{defin}

There are several examples of completely positive kernels, e.g.:
\begin{enumerate}
\item $a(t)= t^{\alpha-1}/\Gamma(\alpha), \quad \alpha\in (0,1],\quad t>0$;
\item $a(t)= e^{-t}, \quad t\ge 0 $.
\end{enumerate}

This is interesting that the class of 
completely positive kernels appears naturally in the theory of
viscoelasticity. Several properties and examples of such kernels
appears in \cite[Section 4.2]{Pr2}.

\begin{defin}\label{dSW3}
Suppose $S(t),~t\geq 0$, is a resolvent for the equation (\ref{eSW1d}). $S(t)$
is called {\tt expo\-nen\-tially bounded} if there are constants
$M\geq 1$ and $\omega\in\mathbb{R}$ such that
$$ ||S(t)|| \leq M\,e^{\omega t}, \mbox{~~for all~~} t\geq 0. $$
\end{defin}

Let us note that in contrary to the semigroup case, not every
resolvent needs to be exponentially bounded even if the kernel
function $a$ belongs to $L^1(\mathbb{R}_+;\mathbb{R})$. 
The class of
kernels providing such  resolvents are $a(t)=
t^{\beta-1}/\Gamma(\beta),~\alpha\in (0,2)$ or the class of
completely monotonic functions.

In this paper the following results concerning convergence of
resolvents for the equation (\ref{eSW1}) in a Hilbert space 
will play the key role. These results have been proved by  A. Karczewska and C. Lizama 
in  \cite{KaLi2}.

\begin{theor} \label{pSW2} (see \ \cite{KaLi2})\\
Let $A$ be the generator of a $C_0$-semigroup in $H$ and suppose
the kernel function $a$ is completely positive. Then the pair $(A,a)$
admits an exponentially bounded resolvent $S(t)$. Moreover, there
exist bounded operators $A_n$ such that $(A_n,a)$ admit resolvent
families $S_n(t)$ satisfying $ ||S_n(t) || \leq Me^{w_0 t}~ (M\geq
1,~w_0\geq 0)$ for all $t\geq 0$ and
\begin{equation} \label{eSW4}
S_n(t)x \to S(t)x \quad \mbox{as} \quad n\to +\infty
\end{equation}
for all $x \in H,\; t\geq 0.$ Additionally, the convergence is
uniform in $t$ on every compact subset of $ \mathbb{R}_+$.
\end{theor}

In the paper the operators 
\begin{equation} \label{eSW5}
A_n := n AR(n,A) = n^2 R(n,A) - nI, \qquad n> w
\end{equation}
are the {\tt Yosida approximation} of $A$.

An analogous result like Theorem \ref{pSW2} holds in other cases.

\begin{theor} \label{pSW2a} (see \ \cite{KaLi2})\\
 Let $A$ be the generator of a strongly continuous cosine family.
 Suppose any of the following:
 \begin{description}
  \item[\hspace{1ex}(i)] $a\in L_\mathrm{loc}^1(\mathbb{R}_+;\mathbb{R})$ is completely
  positive;
  \item[\hspace{0.5ex}(ii)] the kernel fuction $a$ is a creep function with $a_1$
  log-convex;
  \item[(iii)] $a=c\star c$ with some completely positive
  $c\in L_\mathrm{loc}^1(\mathbb{R}_+)$.
 \end{description}
Then the pair $(A,a)$ admits an exponentially bounded resolvent $S(t)$.
Moreover, there exist bounded operators $A_n$ such that $(A_n,a)$
admit resolvent families $S_n(t)$ satisfying $ ||S_n(t) || \leq
Me^{w_0 t}$ $(M\geq 1,~w_0\geq 0)$ for all $t\geq 0,~n\in
\mathbb{N},$ and
$ ~S_n(t)x \to S(t)x$ as $n\to +\infty$
for all $x \in H$, $t\geq 0$. Additionally, the convergence is
uniform in $t$ on every compact subset of $ \mathbb{R}_+$.
\end{theor}

\begin{prop}\label{AScom}
 Let $A, A_n$ and $ S_n( t)$ be given as in Theorem \ref{pSW2}. 
Then $ S_n(t)$ commutes with  the operator $A$, for every
 $n$ sufficiently large and $t\ge 0$.
\end{prop}

\section{The stochastic integral} \label{sSW2}

In this section we recall the construction of the stochastic integral due to van Gaans  \cite{Ga}. 

\begin{defin}
A function  $f:[0,T]\rightarrow L^2(\Omega;H)$ is called {\tt piecewise uni\-form\-ly continuous (PUC)} 
if there are $0=a_0<a_1<\ldots<a_n=T$ such that $f$ is uniformly continuous on  $(a_{k-1},a_k)$ for each $k\in\{1,2,\ldots,n\}$.
\end{defin}

\begin{defin}
A function  $f:[0,\infty)\rightarrow L^2(\Omega;H)$ is called  {\tt  piecewise uni\-form\-ly continuous (PUC)}, if  $f|_{[0,T]}$ is uniformly continuous for all $T>0$.
\end{defin}

\begin{theor} \label{SI1} (see \ \cite{Ga})\\
Assume that $(W_i)_{i=1}^{\infty}$ is a series of independent standard scalar Wiener processes with respect to the filtration $(\mathcal{F}_t)_{t\geq0}$ in $\mathcal{F}$. 
Let $(\Psi_i)_{i=1}^{\infty}$ be a series of piecewise uniformly continuous functions (PUC) 
acting from  $[0,T]$ into $L^2(\Omega;H)$, adapted with respect to the filtration $(\mathcal{F}_t)_{t\geq0}$. Then the following results hold.
\begin{enumerate}
	\item For any  $i\in\mathbb{N}$, the integral $\int\limits_0^T \Psi_i(\tau)dW_i(\tau)$ is well-defined as the limit of Riemann sums $L^2(\Omega;H)$ of the form
	$$\sum\limits_{k=1}^n\Psi_i(\tau_{k-1})(W_i(\tau_k)-W_i(\tau_{k-1})),$$
	where $0=\tau_0<\tau_1<\ldots<\tau_n=T$.\

	\item For each  $i\in\mathbb{N}$, the It\^o isometry holds
	$$\mathbb{E}\left\vert\int\limits_0^T\Psi_i(\tau)dW_i(\tau)\right\vert_H^2 = \int\limits_0^T\mathbb{E}\left\vert\Psi_i(\tau)\right\vert_{L^2(\Omega;H)}^2 d\tau.$$\

	\item For any  $i,j\in\mathbb{N}$, such that $i\neq j$ we have 
	$$\mathbb{E}\left\langle \int\limits_0^T \Psi_i(\tau)dW_i(\tau), \int\limits_0^T \Psi_j(\tau)dW_j(\tau) \right\rangle_H = 0.$$
\end{enumerate}
\end{theor}

\begin{defin}
By  $\mathcal{PUC}([0,T];L^2(\Omega;H))$ we shall denote the space of series  $\Psi=(\Psi_i)_{i=1}^{\infty}$  of piecewise uniformly continuous functions (PUC) acting from $[0,T]$ into $L^2(\Omega;H)$, adapted with respect to the filtration $(\mathcal{F}_t)_{t\geq0}$, such that  
$$\int\limits_0^T\sum\limits_{i=1}^{\infty}\mathbb{E}\vert\Psi_i(\tau)\vert_{L^2(\Omega;H)}^2 d\tau<\infty.$$
\end{defin}

\begin{theor} \label{SI2} (see  \cite{Ga}) \label{st_int-const} \\
Assume that  $W=(W_i)_{i=1}^{\infty}$   is a series of independent standard scalar Wiener processes with respect to the filtration  $(\mathcal{F}_t)_{t\geq0}$ in $\mathcal{F}$. 
Let $\Psi\in \mathcal{PUC} ([0,T];L^2(\Omega;H))$.
Then the integral
$$\int\limits_0^T\Psi(\tau)dW(\tau) := \sum\limits_{i=1}^{\infty}\int\limits_0^T \Psi_i(\tau)dW_i(\tau)$$
exists in  $L^2(\Omega;H)$ and 
$$\mathbb{E}\left\vert\int\limits_0^T\Psi(\tau)dW(\tau)\right\vert_H^2=\int\limits_0^T\mathbb{E}\sum\limits_{i=1}^{\infty}\vert\Psi_i(\tau)\vert_{L^2(\Omega;H)}^2 d\tau.$$
\end{theor}

\section{The main results} \label{sSW3}

We begin this section with definitions of solutions to the equation  (\ref{eSW1}).\

\begin{defin} \label{dSW4}
An $H$-valued predictable process
$X(t),~t\in [0,T]$, is said to be a ~{\tt strong solution}~ to
(\ref{eSW1}), if $X$ has a version such that $P(X(t)\!\in\! D(A))\!=\!1$ 
for almost all $t\in [0,T]$; for any $t\in [0,T]$
\begin{equation} \label{eSW3.1}
\int_0^t |a(t-\tau)AX(\tau)|_H \,d\tau<\infty,\quad P\mbox{--a.s.}
\end{equation}
and for any $t\in [0,T]$ the equation (\ref{eSW1}) holds $P$--a.s.
\end{defin}

Let $A^*$ denote the adjoint of $A$ with a dense domain
$D(A^*)\subset H$ and the graph norm $|\cdot |_{D(A^*)}$.

\begin{defin} \label{dSW5}
An $H$-valued predictable process $X(t),~t\in
[0,T]$, is said to be a {\tt weak solution} to (\ref{eSW1}), if
$P(\int_0^t|a(t-\tau)X(\tau)|_H d\tau<\infty)=1$ and if for all
$\xi\in D(A^*)$ and all $t\in [0,T]$ the following equation holds
\begin{eqnarray*}
\langle X(t),\xi\rangle_H = \langle X_0,\xi\rangle_H &+& \left\langle
\int_0^t a(t-\tau)X(\tau)\,d\tau, A^*\xi\right\rangle_H \\ &+& \left\langle \sum\limits_{i=1}^{\infty}\int\limits_0^t \Psi_i(\tau)dW_i(\tau),\xi\right\rangle_H, \quad P\mbox{--a.s.}
\end{eqnarray*}
\end{defin}

As we have already written, in the paper we assume that (\ref{eSW1d}) admits a resolvent
family $S(t)$, $t \geq 0.$ So, we can introduce the following idea.

\begin{defin} \label{dSW6}
An $H$-valued predictable process $X(t),~t\in
[0,T]$, is said to be a {\tt mild solution} to the stochastic
Volterra equation (\ref{eSW1}), if
\begin{equation}\label{edf6}
 \int\limits_0^t\mathbb{E}\sum\limits_{i=1}^{\infty}\vert S(t-\tau)\Psi_i(\tau)\vert_{L^2(\Omega;H)}^2 d\tau <\infty\ \ {for}\ \ t\leq T
 \end{equation}
 and, for arbitrary $t\in[0,T]$,
\begin{equation}\label{eSW9}
 X(t) = S(t)X_0 + \sum\limits_{i=1}^{\infty}\int\limits_0^t S(t-\tau)\Psi_i(\tau)dW_i(\tau),\quad P\mbox{--a.s.}\end{equation}
where $S(t), t\ge 0$, is the resolvent for the equation (\ref{eSW1d}).
\end{defin}

In the paper we will use the following generalization of the well-known result, e.g., compare Proposition 4.15 in \cite{DaPrZa}.

\begin{prop} \label{prop314}
Assume that  $\Phi_i(t) \in D(A)$, $P\mbox{--a.s.}$ for $i\in\mathbb{N}$ and $t\in[0,T]$. If
\begin{equation} \label{pro314eq1}
\int\limits_0^T\sum\limits_{i=1}^{\infty}\mathbb{E}\vert\Phi_i(t)\vert_{L^2(\Omega;H)}^2 dt<\infty,\ \ \ \int\limits_0^T\sum\limits_{i=1}^{\infty}\mathbb{E}\vert A\Phi_i(t)\vert_{L^2(\Omega;H)}^2 dt<\infty,
\end{equation}
then
\begin{equation}\label{pro314eq2}
 P\left(\sum\limits_{i=1}^{\infty}\int\limits_0^T \Phi_i(t)dW_i(t) \in D(A)\right)=1
\end{equation}
and
$$ A\sum\limits_{i=1}^{\infty}\int\limits_0^T \Phi_i(t)dW_i(t) = \sum\limits_{i=1}^{\infty}\int\limits_0^T A\Phi_i(t)dW_i(t),\quad P\mbox{--a.s.}$$
\end{prop}

\bgproof
From the assumption we have that  $\Phi_i(t) \in D(A)$, $P-a.s.$ for any $i\in\mathbb{N}$ and $t\in[0,T]$. 
So, the conditions  (\ref{pro314eq1}) can be written as 
\begin{equation} \label{pro314eq1a}
\int\limits_0^T\sum\limits_{i=1}^{\infty}\mathbb{E}\vert\Phi_i(t)\vert_{L^2(\Omega;D(A))}^2 dt<\infty,\ \ \ \int\limits_0^T\sum\limits_{i=1}^{\infty}\mathbb{E}\vert A\Phi_i(t)\vert_{L^2(\Omega;D(A))}^2 dt<\infty.
\end{equation}
Because $(D(A),\vert\cdot\vert_{D(A)})$ is a Hilbert space, then the integral 
$$\int\limits_0^T\Phi(t)dW(t) := \sum\limits_{i=1}^{\infty}\int\limits_0^T \Phi_i(t)dW_i(t)$$
exists in  $L^2(\Omega;D(A))$  by Theorem \ref{st_int-const}.

Denote by  $(\tilde{s}_n)_{n=1}^{\infty}$ division of the interval  $[0,T]$, that $0=s_0<s_1<\ldots<s_n=T$. From the definition of the integral and closedness of the operator $A$ 
we have
\begin{eqnarray*}
\sum\limits_{i=1}^{\infty}\int\limits_0^T A\Phi_i(t)dW_i(t) & = & \sum\limits_{i=1}^{\infty}\lim\limits_{n\rightarrow\infty} A\Phi_i(s_n)(W_i(s_n)-W_i(s_{n-1})) \\
 &=& \lim\limits_{m\rightarrow\infty}\sum\limits_{i=1}^{m}\lim\limits_{n\rightarrow\infty} A\Phi_i(s_n)(W_i(s_n)-W_i(s_{n-1})) \\
 &=& A\lim\limits_{m\rightarrow\infty}\sum\limits_{i=1}^{m}\lim\limits_{n\rightarrow\infty} \Phi_i(s_n)(W_i(s_n)-W_i(s_{n-1}))  \\
 &=& A\lim\limits_{m\rightarrow\infty}\sum\limits_{i=1}^{m}\int\limits_0^T \Phi_i(t)dW_i(t) \\ &=&
 A\sum\limits_{i=1}^{\infty}\int\limits_0^T \Phi_i(t)dW_i(t),\quad P\mbox{--a.s.}
\end{eqnarray*}
\edproof

 We introduce the {\tt stochastic convolution}
\begin{equation} \label{eSW18a}
W^{\Psi}(t):=\sum\limits_{i=1}^{\infty}\int\limits_0^t S(t-\tau)\Psi_i(\tau)dW_i(\tau),\ t\in[0,T],
\end{equation}
where $\Psi$ and the resolvent operators $S(t),~t\geq 0$, are the same
as above.

Let us formulate some auxiliary results concerning the convolution
$W^\Psi$.

\begin{prop}\label{pr2proper}
For arbitrary process  $\Psi\in \mathcal{PUC}([0,T];L^2(\Omega;H))$,
the process $W^\Psi(t)$, $t\geq 0$, given by (\ref{eSW18a}) has a
predictable version.
\end{prop}

\begin{prop}\label{pr3proper}
Assume that $\Psi\in \mathcal{PUC}([0,T];L^2(\Omega;H))$. Then the process
$W^\Psi(t)$, $t\geq 0$, defined by (\ref{eSW18a}) has square
integrable trajectories.
\end{prop}

For  the idea of proofs of Propositions \ref{pr2proper} and
\ref{pr3proper} we refer to \cite{Ka1} or  \cite{Ka2}.

In some cases weak solutions of equation (\ref{eSW1}) coincides
with mild solutions of (\ref{eSW1}), see e.g.\ \cite{Ka1} or \ \cite{Ka2}  . In
consequence, having results for the convolution  (\ref{eSW18a}) we obtain results for weak solutions.

\begin{prop} \label{pr4proper}
Let $a\in BV(\mathbb{R}_+;\mathbb{R})$ and suppose that (\ref{eSW1d}) admits
a resolvent family $S \in C^1(0,\infty; L(H)).$ Let $X$ be a
predictable process with integrable trajectories. Assume that $X$
has a version such that $P(X(t)\in D(A))=1$ for almost all $t\in
[0,T]$ and (\ref{edf6}) holds. If for any $t\in [0,T]$ and $\xi\in
D(A^*)$
\begin{eqnarray}\label{deq9}
 \langle X(t),\xi\rangle_H = \langle X_0,\xi\rangle_H &+&
 \int_0^t \langle a(t-\tau)X(\tau),A^*\xi\rangle_H d\tau \\
 &+& \left\langle\sum\limits_{i=1}^{\infty}\int\limits_0^t \Psi_i(\tau)dW_i(\tau),\xi\right\rangle_H, 
  ~P-a.s., \nonumber
\end{eqnarray}
then
\begin{equation}\label{deq9a}
X(t) = S(t)X_0 +
\sum\limits_{i=1}^{\infty}\int\limits_0^t S(t-\tau)\Psi_i(\tau)dW_i(\tau),\ t\in[0,T].
\end{equation}
\end{prop}

\begin{rem}
{\em 
If  (\ref{eSW1}) is parabolic and the kernel $a$ is
3-monotone, understood in the sense defined by Pr\"uss
\cite[Section 3]{Pr2}, then $S \in C^1(0,\infty; L(H))$ and $a\in
BV( \mathbb{R}_+;\mathbb{R})$, respectively. }
\end{rem}

\begin{prop} \label{pr5prop}
Assume that $A$, the kernel function and the operators 
$S(t)$, $t\ge 0$, are as previously.
If $\Psi\in \mathcal{PUC}([0,T];L^2(\Omega;H)),$  then the stochastic convolution
$W^\Psi$ fulfills the equation (\ref{deq9}) with $X_0\equiv 0$.    
\end{prop}

Hence, we are able to conclude the following result.\\

\begin{cor} \label{c1prop}
Assume that $A$ is a linear bounded operator in $H$, \linebreak
$a\in L_\mathrm{loc}^1(\mathbb{R}_+;\mathbb{R})$ and $S(t)$, $t\ge 0$, 
are resolvent operators for the equation~(\ref{eSW1d}).
If $\Psi\in \mathcal{PUC}([0,T];L^2(\Omega;H)),$ then
\begin{equation} \label{deq16}
 W^\Psi(t) = \int_0^t a(t-\tau)AW^\Psi(\tau)d\tau
 +\sum\limits_{i=1}^{\infty}\int\limits_0^t \Psi_i(\tau)dW_i(\tau)\;,
\quad\quad t\geq 0\;.
\end{equation}
\end{cor}

\begin{rem} {\em The formula (\ref{deq16}) says that the convolution $W^\Psi$
is a strong solution to (\ref{eSW1}) with $X_0\equiv 0$ if the
operator $A$ is bounded.}
\end{rem}

Here we provide sufficient conditions under which the
stochastic convolution $W^\Psi (t)$, $t\ge 0$, defined by
(\ref{eSW18a}) is a strong solution to the equation (\ref{eSW1}).\\

\begin{theor} \label{pSW5}
Let $A$ be a closed linear unbounded operator with the dense
domain $D(A)$ equipped with the graph norm $|\cdot|_{D(A)}$.
Suppose that assumptions of Theorem \ref{pSW2} or Theorem
\ref{pSW2a} hold. If  $\Psi$ and $A\Psi$ belong to  $\mathcal{PUC}([0,T];L^2(\Omega;H))$  and for every  $i\in\mathbb{N}$, $\Psi_i(t)\in D(A)$, $P\mbox{--a.s.}$, then the following equality holds
\begin{equation}\label{lemma_equation}
W^{\Psi}(t) = \int\limits_0^t a(t-\tau)AW^{\Psi}(\tau)d\tau + \sum\limits_{i=1}^{\infty}\int\limits_0^t \Psi_i(\tau)dW_i(\tau),\quad P\mbox{--a.s.}
\end{equation}
\end{theor}

For reader's convenience we postponed the proof of Theorem \ref{pSW5} 
to section \ref{sSW5}.

\begin{theor} \label{coSW6}
Suppose that assumptions of Theorem \ref{pSW2} or Theorem
\ref{pSW2a} hold. Then the equation (\ref{eSW1}) has a strong
solution. Precisely, the convolution $W^\Psi$ defined by
(\ref{eSW18a}) is the strong solution to~(\ref{eSW1}) with
$X_0\equiv 0$.
\end{theor}
\bgproof
In order to prove Theorem \ref{coSW6}, we have to show
only the condition (\ref{eSW3.1}). Let us note that the
convolution $W^\Psi$ has integrable trajectories. Because the
closed unbounded linear 
operator $A$ becomes bounded on ($D(A),|\cdot|_{D(A)}$), see
\cite[Chapter 5]{We}, we obtain that $AW^\Psi (\cdot )\in
L^1([0,T];H)$, $~P$--a.s. Hence, properties of convolution provide
integrability of the function $a(T-\tau)AW^\Psi (\tau)$ with
respect to $\tau$, what finishes the proof. 
\edproof

\section{Regularity of stochastic convolution} \label{sSW4}

In this section we prove thet the stochastic convolution given by (\ref{eSW18a}) has continuous trajectories. In order to do this, we rewrite the convolution to more convenient form (see the formula (\ref{e41}) below). Then, using the obtained formula we join it with an appropriate Cauchy problem. Next, we use well-known results for some Cauchy problems considered, e.g.\ in \cite{Pa} and we adapt them for our purposes.

\begin{theor}\label{th_reg_1}
Assume that the operator $A$ is the generator of a strongly continuous semigroup $\tilde{T}(t)$, $t\in[0,T]$. Let the kernel function $a(t)$, $t\in [0,T]$, be completely positive, $\dot{a} \in L^1_{loc}([0,T];\mathbb{R})$
and $\Psi\in \mathcal{PUC}([0,T];L^2(\Omega;H))$. Then the following equation holds
\begin{eqnarray} \label{e41}
W^{\Psi}(t) &=& cA\int\limits_0^t \tilde{T}(t-\tau)\left[\int\limits_0^{\tau}\dot{a}(\tau-\sigma)W^{\Psi}(\sigma)d\sigma + c\int\limits_0^{\tau}\Psi(s)dW(s)\right]d\tau \nonumber \\ &&+ \int\limits_0^{t}\Psi(\tau)dW(\tau),
\end{eqnarray}
where $t\in[0,T]$ and $c=a(0)$ is a constant.
\end{theor}

\bgproof
Because the equation (\ref{deq16}) holds for any bounded operator, it holds also for the Yosida approximation $A_n$ of the operator $A$, that is
\begin{equation}\label{reg_proof_1}
W^{\Psi}_n(t) = \int\limits_0^t a(t-\tau)A_nW^{\Psi}_n(\tau)d\tau + \sum\limits_{i=1}^{\infty}\int\limits_0^t \Psi_i(\tau)dW_i(\tau),
\end{equation}
where
$$
W^{\Psi}_n(t) := \sum\limits_{i=1}^{\infty}\int\limits_0^t S_n(t-\tau)\Psi_i(\tau)dW_i(\tau).
$$
Let us define 
\begin{equation}\label{reg_proof_2}
Z_n(t):=\int\limits_0^t a(t-\tau)W^{\Psi}_n(\tau)d\tau, \quad t\in[0,T]\,;
\end{equation}
 using the Leibniz rule we obtain
\begin{equation}\label{reg_proof_3}
Z_n'(t)=\int\limits_0^t \dot{a}(t-\tau)W^{\Psi}_n(\tau)d\tau + a(0)W^{\Psi}_n(t), \quad t\in[0,T]\,.
\end{equation}
From 
 (\ref{reg_proof_1}) and (\ref{reg_proof_2}) we have
$$
W^{\Psi}_n(t)=A_nZ_n(t)+\int\limits_0^{t}\Psi(\tau)dW(\tau), \ \ t\in[0,T].
$$
If $a(0)\neq 0$, then from (\ref{reg_proof_3}) we receive
$$
W^{\Psi}_n(t)=\frac{1}{a(0)}\left[Z_n'(t)- \int\limits_0^t \dot{a}(t-\tau)W^{\Psi}_n(\tau)d\tau \right],
$$
that is
$$
Z_n'(t)- \int\limits_0^t \dot{a}(t-\tau)W^{\Psi}_n(\tau)d\tau = a(0)[A_nZ_n(t)+\int\limits_0^{t}\Psi(\tau)dW(\tau)], \quad t\in[0,T].
$$
This means that
$$
Z_n'(t) = a(0)A_nZ_n(t)+\int\limits_0^t \dot{a}(t-\tau)W^{\Psi}_n(\tau)d\tau+ a(0)\int\limits_0^{t}\Psi(\tau)dW(\tau), \ \ t\in[0,T].
$$
Let us denote
$$
\tilde{W}^{\Psi}_n(t) := \int\limits_0^t \dot{a}(t-\tau)W^{\Psi}_n(\tau)d\tau, \ \ t\in[0,T].
$$
Then
$$
Z_n'(t) = cA_nZ_n(t)+ [\tilde{W}^{\Psi}_n(t)+ c\int\limits_0^{t}\Psi(\tau)dW(\tau)],\quad \textrm{where}\quad c=a(0).
$$
Because $Z_n(0)=0$, then
$$
Z_n(t)=\int\limits_0^t e^{c(t-\tau)A_n}\left[\tilde{W}^{\Psi}_n(\tau)+ c\int\limits_0^{\tau}\Psi(s)dW(s)\right]d\tau, \ \ t\in[0,T].
$$
From (\ref{reg_proof_1}) one has
$$
W^{\Psi}_n(t)=AJ_nZ_n(t)+\int\limits_0^{t}\Psi(\tau)dW(\tau), \ \ t\in[0,T],
$$
where $J_n:=nR(n,A)$. 

Therefore
\begin{eqnarray*}
W^{\Psi}_n(t) &=& AJ_n\int\limits_0^t e^{c(t-\tau)A_n}\left[\tilde{W}^{\Psi}_n(\tau)+ c\int\limits_0^{\tau}\Psi(s)dW(s)\right]d\tau \\ &&
+\int\limits_0^{t}\Psi(\tau)dW(\tau), \quad t\in[0,T].
\end{eqnarray*}
By Theorem \ref{pSW2}, properties of the Yosida approximation $A_n$ of the operator $A$ and the Lebesgue dominated convergence theorem one obtains 
\begin{eqnarray*}
\lim\limits_{n\rightarrow\infty}J_n x & = & x, \textrm{ for arbitrary } x\in H, \\
\lim\limits_{n\rightarrow\infty}A_n x &  = &  Ax,  \textrm{ for arbitrary } x\in D(A),\\
\lim\limits_{n\rightarrow\infty}e^{tA_n}x  & = &  \tilde{T}(y)x,  \textrm{ for arbitrary } x\in H,\\
\textrm{and} && \lim\limits_{n\rightarrow\infty}\sup\limits_{t\in[0,T]}\mathbb{E}\left\vert W_n^{\Psi}(t)-W^{\Psi}(t)\right\vert_H^2=0.
\end{eqnarray*}
Because the operator $A$ is closed we conclude that
$$
\int\limits_0^t \tilde{T}(t-\tau)\left[\tilde{W}^{\Psi}_n(\tau)+ c\int\limits_0^{\tau}\Psi(s)dW(s)\right]d\tau \in D(A).
$$
When $n\rightarrow\infty$, we get
\begin{eqnarray*}
W^{\Psi}(t) &=& cA\int\limits_0^t \tilde{T}(t-\tau)\left[\tilde{W}^{\Psi}(\tau)+ c\int\limits_0^{\tau}\Psi(s)dW(s)\right]d\tau \\ &&
+\int\limits_0^{t}\Psi(\tau)dW(\tau), \quad t\in[0,T],
\end{eqnarray*}
where
$$
\tilde{W}^{\Psi}(\tau) = \int\limits_0^\tau \dot{a}(\tau-\sigma)W^{\Psi}(\sigma)d\sigma, \quad \tau\in[0,T].
$$
\edproof

\begin{cor} \label{c42}
Let assumptions of Theorem  \ref{th_reg_1} hold and let $a(0)=1$. Then, for $t\in[0,T]$,
\begin{equation}
W^{\Psi}(t)=AY(t)+\int\limits_0^{t}\Psi(\tau)dW(\tau),
\end{equation}
where
\begin{equation}\label{ydef}
Y(t):=\int\limits_0^t \tilde{T}(t-\tau)\left[\tilde{W}^{\Psi}(\tau)+ \int\limits_0^{\tau}\Psi(s)dW(s)\right]d\tau,
\end{equation}
and
$$
\tilde{W}^{\Psi}(\tau) = \int\limits_0^\tau \dot{a}(\tau-\sigma)W^{\Psi}(\sigma)d\sigma, \ \ \tau\in[0,T].
$$
Moreover, $Y\in C^1([0,T];D(A))$, $P$--a.s. and
\begin{equation} \label{e47}
\frac{dY(t)}{dt}=AY(t)+\left[\tilde{W}^{\Psi}(t)+ \int\limits_0^{t}\Psi(s)dW(s)\right], \ \ t\in[0,T].
\end{equation}
\end{cor}
\bgproof
From the definition (\ref{ydef}) of the process $Y$ and properties of the convolution, $Y\in C^1([0,T];D(A))$, $P$--a.s. Using the Leibniz rule and semigroup property we have
\begin{eqnarray*}
\frac{dY(t)}{dt} & = & \int\limits_0^t \frac{dT(t-\tau)}{dt}\left[\tilde{W}^{\Psi}(\tau)+ \int\limits_0^{\tau}\Psi(s)dW(s)\right]d\tau \\ && + T(0)\left[\tilde{W}^{\Psi}(t)+ \int\limits_0^{t}\Psi(s)dW(s)\right]\\
 & = & A\int\limits_0^t T(t-\tau)\left[\tilde{W}^{\Psi}(\tau)+ \int\limits_0^{\tau}\Psi(s)dW(s)\right]d\tau \\ && + \left[\tilde{W}^{\Psi}(t)+ \int\limits_0^{t}\Psi(s)dW(s)\right]\\
 & = & AY(t) + \left[\tilde{W}^{\Psi}(t)+ \int\limits_0^{t}\Psi(s)dW(s)\right], \quad t\in[0,T]\,.
\end{eqnarray*}
\edproof

This is worth to emphasize that Theorem \ref{th_reg_1} and Corrolary \ref{c42} can be reformulated under different assumptions on the kernel function $a(t), t\ge 0$, and the operator $A$. If $a(t), t\ge 0$, and $A$ satisfy the assumptions of Theorem~\ref{pSW2a}, we obtain the same results as above.

We finish this section with the following regularity result.

\begin{theor}
Let assumptions of the Theorem \ref{th_reg_1} are fulfilled and let the semigroup generated by the operator $A$ is the analytical one. Then the stochastic convolution $W^{\Psi}(t)$, $t\geq 0$, given by (\ref{eSW18a}) has continuous  trajectories.
\end{theor}
\bgproof
Let us note that the process $Y(t), t\ge 0$, and the stochastic integral $\int_0^t \psi(\tau)dW(\tau)$ have continuous trajectories. Hence, using the regularity results for nonhomogeneous Cauchy problem from \cite{Pa} to the formula (\ref{e47}), we conclude continuity of trajectories of the stochastic convolution (\ref{eSW18a}).
\edproof

\section{Proof of Theorem \ref{pSW5}} \label{sSW5}

For the reader's convenience we write the proof with details.

As we have already written, 
because the formula (\ref{deq16})  holds for any bounded operator, then it holds for the Yosida 
approximation  $A_n$ of the operator $A$, too. So, we have 
$$
W^{\Psi}_n(t) = \int\limits_0^t a(t-\tau)A_nW^{\Psi}_n(\tau)d\tau + \sum\limits_{i=1}^{\infty}\int\limits_0^t \Psi_i(\tau)dW_i(\tau),
$$
where
$$
W^{\Psi}_n(t) := \sum\limits_{i=1}^{\infty}\int\limits_0^t S_n(t-\tau)\Psi_i(\tau)dW_i(\tau)
$$
and
$$
A_nW^{\Psi}_n(t) := A_n\sum\limits_{i=1}^{\infty}\int\limits_0^t S_n(t-\tau)\Psi_i(\tau)dW_i(\tau).
$$
By assumption  $\Psi \in \mathcal{PUC}([0,T];L^2(\Omega;H))$. Because the operators  $S_n(t)$ are deterministic and bounded for any  $t\in [0,T]$, $n\in\mathbb{N}$, then $S_n(t-\cdot)\Psi(\cdot)$ belong to  $\mathcal{PUC}([0,T];L^2(\Omega;H))$, too. Hence, the difference
\begin{equation}\label{difference}
\Phi_{n}(t-\cdot) := S_n(t-\cdot)\Psi(\cdot)-S(t-\cdot)\Psi(\cdot)
\end{equation}
belongs to $\mathcal{PUC}([0,T];L^2(\Omega;H))$ for every $t\in [0,T]$ and $n\in\mathbb{N}$, too. 
This means that 
\begin{equation}\label{int_diff}
\int\limits_0^t\sum\limits_{i=1}^{\infty}\mathbb{E}\vert\Phi_{n,i}(t-\tau)\vert_{L^2(\Omega;H)}^2 d\tau<\infty \textrm{ \ \ for } t\in [0,T].
\end{equation}
From the definiton of the stochastic integral introduced in Theorem \ref{st_int-const}, for  $t\in[0,T]$ we have 
$$\mathbb{E}\left\vert\int\limits_0^t\Phi_n(t-\tau)dW(\tau)\right\vert_H^2=\int\limits_0^t\mathbb{E}\sum\limits_{i=1}^{\infty}\vert\Phi_{n,i}(t-\tau)\vert_{L^2(\Omega;H)}^2 d\tau < \infty.$$

By the approximation theorems, that is Theorem \ref{pSW2} or Theorem \ref{pSW2a},
the convergence of the resolvent families is uniform with respect to $t$ on every closed intervals, particularly on. Then we have 
\begin{equation}\label{res_convergence}
\int\limits_0^T\mathbb{E}\sum\limits_{i=1}^{\infty}\vert [S_n(T-\tau)-S(T-\tau)]\Psi_i(\tau)\vert_{L^2(\Omega;H)}^2 d\tau \rightarrow 0 \quad \textrm{for} \quad n\rightarrow +\infty.
\end{equation}

Summing up the above considerations, we obtain 
\begin{eqnarray*}
\sup\limits_{t\in[0,T]} &\mathbb{E}&\left\vert\int\limits_0^t\Phi_n(t-\tau)dW(\tau)\right\vert_H^2  \equiv \\
&\equiv& \sup\limits_{t\in[0,T]} \mathbb{E}\left\vert\int\limits_0^t \left[S_n(t-\tau)-S(t-\tau)\right]\Psi(\tau)dW(\tau)\right\vert_H^2 \\
& = & \sup\limits_{t\in[0,T]} \mathbb{E}\left\vert\sum\limits_{i=1}^{\infty}\int\limits_0^t[S_n(t-\tau)-S(t-\tau)]\Psi_i(\tau)dW_i(\tau)\right\vert_H^2 \\
& \leq & \mathbb{E}\left\vert\sum\limits_{i=1}^{\infty}\int\limits_0^T[S_n(T-\tau)-S(T-\tau)]\Psi_i(\tau)dW_i(\tau)\right\vert_H^2 \\
& = &\int\limits_0^T\mathbb{E}\sum\limits_{i=1}^{\infty}\vert [S_n(T-\tau)-S(T-\tau)]\Psi_i(\tau)\vert_{L^2(\Omega;H)}^2 d\tau \rightarrow 0
\end{eqnarray*}
as $n\rightarrow +\infty$.

By the dominated convergence theorem we can obtain 
\begin{equation}\label{convergence1}
\lim\limits_{n\rightarrow\infty}\sup\limits_{t\in[0,T]}\mathbb{E}\left\vert W_n^{\Psi}(t)-W^{\Psi}(t)\right\vert_H^2=0.
\end{equation}

By assumption, $\Psi_i(t)\in D(A)$, for $i\in\mathbb{N}$, $P--a.s$.
Because $S(t)(D(A))\subset D(A)$, so $S(t-\tau)\Psi_i(\tau)(H)\subset D(A)$, $P$--a.s., for any $\tau\in [0,t]$, $t\geq 0$, $i\in\mathbb{N}$. 
Then, by Proposition \ref{prop314}, $P(W^{\Psi}(t)\in D(A))=1$.

For any  $n\in\mathbb{N}$, $t\geq 0$, we can write
$$ \vert A_n W_n^{\Psi}(t)-AW^{\Psi}(t)\vert_H \leq N_{n,1}(t) + N_{n,2}(t),$$
where
$$N_{n,1}(t):=  \vert A_n W_n^{\Psi}(t)-A_n W^{\Psi}(t)\vert_H \qquad \mbox{and}$$
$$N_{n,2}(t):=  \vert A_n W^{\Psi}(t)-AW^{\Psi}(t)\vert_H=\vert (A_n-A)W^{\Psi}(t)\vert_H. $$
So
\begin{eqnarray}\label{eq329}
\vert A_n W_n^{\Psi}(t)-AW^{\Psi}(t)\vert^2_H &\le & N^2_{n,1}(t) + 2N_{n,1}(t)N_{n,2}(t)+ N^2_{n,2}(t) \nonumber\\ 
&\le & 3[N^2_{n,1}(t)+N^2_{n,2}(t)].
\end{eqnarray}

We shall study the term $N_{n,1}(t)$ first. Because the operator $A$ generates a semigroup, we can use the following property of the Yosida approximation:
\begin{equation}\label{yosida_properties1}
 A_n x=J_n Ax \textrm{ for any } x\in D(A),\ \sup\limits_n\Vert J_n\Vert <\infty, 
\end{equation}
where $A_n x=nAR(n,A)x=AJ_n x$ for any $x\in H$, $J_n:=nR(n,A)$.

 Moreover
\begin{eqnarray}\label{yosida_convergence}
\lim\limits_{n\rightarrow\infty}J_n x = x & \textrm{ for every } x\in H, \\
\lim\limits_{n\rightarrow\infty}A_n x = Ax & \textrm{ for every } x\in D(A).\nonumber
\end{eqnarray}
By Proposition \ref{AScom} for any big enough $n$ and any  $x\in D(A)$ we have \linebreak $AS_n(t)x=S_n(t)Ax$.

Next, by Proposition \ref{prop314} and closedness of the operator $A$ we can write
\begin{eqnarray*}
A_n W_n^{\Psi}(t) & \equiv & A_n \sum\limits_{i=1}^{\infty}\int\limits_0^t S_n(t-\tau)\Psi_i(\tau)dW_i(\tau) \\
 & = & J_n \sum\limits_{i=1}^{\infty}\int\limits_0^t AS_n(t-\tau)\Psi_i(\tau)dW_i(\tau) \\
& = & J_n \left[\sum\limits_{i=1}^{\infty}\int\limits_0^t S_n(t-\tau)A\Psi_i(\tau)dW_i(\tau)\right].
\end{eqnarray*}

Analogously, we have 
$$ A_n W^{\Psi}(t) = J_n\left[\sum\limits_{i=1}^{\infty}\int\limits_0^t S(t-\tau)A\Psi_i(\tau)dW_i(\tau)\right]. $$
Using  (\ref{yosida_properties1}), we receive 
\begin{eqnarray*}
N_{n,1}(t) & = & \left\vert J_n \sum\limits_{i=1}^{\infty}\int\limits_0^t [S_n(t-\tau)-S(t-\tau)]A\Psi_i(\tau)dW_i(\tau)\right\vert_H \\
 & \leq & \left\vert \sum\limits_{i=1}^{\infty}\int\limits_0^t [S_n(t-\tau)-S(t-\tau)]A\Psi_i(\tau)dW_i(\tau)\right\vert_H.
\end{eqnarray*}

From assumption $A\Psi \in \mathcal{PUC}([0,T];L^2(\Omega;H))$, so the term  $[S_n(t-\cdot)-S(t-\cdot)]A\Psi(\cdot)$ may be treated like the difference $\Phi_n$ defined in  (\ref{difference}).

Then, using  (\ref{yosida_properties1}) and (\ref{convergence1}), for $N^2_{n,1}(t)$ we have
$$ \lim\limits_{n\rightarrow\infty} \sup\limits_{t\in[0,T]}\mathbb{E}(N^2_{n,1}(t)) = 0. $$

For $N^2_{n,2}(t)$ the proof of the convergence (\ref{convergence1}) is the same.
\begin{eqnarray*}
N_{n,2}(t) \!& \!= \!&\! \vert A_n W^{\Psi}(t)-AW^{\Psi}(t)\vert_H \\
 \!& \!\equiv \!&\! \left\vert A_n \sum\limits_{i=1}^{\infty}\int\limits_0^t S(t-\tau)\Psi_i(\tau)dW_i(\tau) - A \sum\limits_{i=1}^{\infty}\int\limits_0^t S(t-\tau)\Psi_i(\tau)dW_i(\tau)\right\vert_H \\
 \!& \!= \!&\! \left\vert \sum\limits_{i=1}^{\infty}\int\limits_0^t [A_n-A]S(t-\tau)\Psi_i(\tau)dW_i(\tau) \right\vert_H.
\end{eqnarray*}

By assumption, $\Psi$ and $A\Psi$ belong to $\mathcal{PUC}([0,T];L^2(\Omega;H))$. Because $A_n$ and $S(t)$, $t\geq 0$ are bounded, so $A_nS(t-\cdot)\Psi(\cdot)\in \mathcal{PUC}([0,T];L^2(\Omega;H))$. 

Analogously,  $AS(t-\cdot)\Psi(\cdot)=S(t-\cdot)A\Psi(\cdot)\in \mathcal{PUC}([0,T];L^2(\Omega;H))$.

We may deduce that  $[A_n-A]S(t-\cdot)\Psi(\cdot)\in \mathcal{PUC}([0,T];L^2(\Omega;H))$, what means that this term can be treated like the difference $\Phi_n$ defined by  (\ref{difference}). 
Hence, for any $t\in[0,T]$, we may write
\begin{eqnarray*}
\mathbb{E}\left(N^2_{n,2}(t)\right) & = & \mathbb{E}\left\vert \sum\limits_{i=1}^{\infty}\int\limits_0^t [A_n-A]S(t-\tau)\Psi_i(\tau)dW_i(\tau) \right\vert^2_H \\
 & \leq & \int\limits_0^T\mathbb{E}\sum\limits_{i=1}^{\infty}\vert[A_n-A]S(T-\tau)\Psi_i(\tau)\vert_{L^2(\Omega;H)}^2 d\tau < \infty.
\end{eqnarray*}

Using the convergence  (\ref{yosida_convergence}), we have
$$ \int\limits_0^T\mathbb{E}\sum\limits_{i=1}^{\infty}\vert[A_n-A]S(T-\tau)\Psi_i(\tau)\vert_{L^2(\Omega;H)}^2 d\tau \rightarrow 0 \quad \mbox{for} \quad n\rightarrow +\infty. $$
Then 
$$ \lim\limits_{n\rightarrow\infty} \sup\limits_{t\in[0,T]}\mathbb{E}(N^2_{n,2}(t)) = 0. $$

From the convergences obtained for  $N^2_{n,1}(t)$ and $N^2_{n,2}(t)$ we can conclude, that
\begin{equation}\label{convergence2}
\lim\limits_{n\rightarrow\infty} \sup\limits_{t\in[0,T]}\mathbb{E}\vert A_n W_n^{\Psi}(t)-AW^{\Psi}(t)\vert^2_H = 0.
\end{equation}

By  (\ref{convergence1}) and  (\ref{convergence2}) we see, that the equality (\ref{lemma_equation}) holds.
\hfill $\square$

\end{document}